\documentclass[aos,MSNbibl,seceqn,citesort,dvips]{arximspdf}
\usepackage{mathrsfs}

\doi{10.1214/10-AOS872}
\volume{39}
\issue{3}
\pubyear{2011}
\firstpage{1471}
\lastpage{1495}

\makeatletter
\newtheorem{teo}{Theorem}[section]
\newtheorem{lem}{Lemma}[section]
\newtheorem{prop}{Proposition}[section]
\newproclaim{cond}{Condition}
\newcommand{\cal}{\mathcal}
\newcommand{\idotsint}{\int\cdots\int}
\newcommand{\iint}{\int\!\!\int}
\makeatother

\begin{document}
\begin{frontmatter}

\title{Asymptotic equivalence of functional linear regression and a
white noise inverse problem}
\runtitle{Asymptotic equivalence}

\begin{aug}
\author[A]{\fnms{Alexander} \snm{Meister}\corref{}\ead[label=e1]{alexander.meister@uni-rostock.de}}
\runauthor{A. Meister}
\affiliation{Universit\"{a}t Rostock}
\address[A]{Institut f\"{u}r Mathematik \\
Universit\"{a}t Rostock \\
Ulmenstrasse 69 \\ D-18051 Rostock\\ Germany\\
\printead{e1}} %adresu isvedimo komanda gale!
\end{aug}

% HISTORY:
\received{\smonth{6} \syear{2010}}

% ABSTRACT
%
\begin{abstract}
We consider the statistical experiment of functional linear regression
(FLR). Furthermore, we introduce a white noise model where one observes
an It\^{o} process, which contains the covariance operator of the
corresponding FLR model in its construction. We prove asymptotic
equivalence of FLR and this white noise model in LeCam's sense under
known design distribution. Moreover, we show equivalence of FLR and an
empirical version of the white noise model for finite sample sizes. As
an application, we derive sharp minimax constants in the FLR model
which are still valid in the case of unknown design distribution.
\end{abstract}

% KEYWORDS
%

\end{frontmatter}

%s1 ###
\section{Introduction}
\label{S1}

We consider the statistical problem of functional linear regression
(FLR). In its standard version, one observes the data $(\mathbf{X},\mathbf{Y})$ where $\mathbf{X} = (X_1,\ldots ,X_n)^T$ are i.i.d. random variables
taking their values in $C([0,1])$, that is, the set consisting of all
continuous functions on the interval $[0,1]$, and $\mathbf{Y} = (Y_1,\ldots
,Y_n)^T$ with
%
%e1.1 ###
\begin{equation} \label{eq:1.1}
Y_j   =   \langle X_j , \theta\rangle  +   \varepsilon_j , \qquad
j=1,\ldots ,n ,
\end{equation}
 where $\langle\cdot, \cdot\rangle$ denotes the
$L_2([0,1])$-inner product throughout this work. The i.i.d. error
variables $\varepsilon_j$ are assumed to be centered and normally
distributed with the variance $\sigma^2$. Moreover, all $X_1,\varepsilon
_1,\ldots ,X_n,\varepsilon_n$ are independent. The goal is to estimate
the regression function $\theta\in\Theta\subseteq L_2([0,1])$. In
general, we allow for such a structure of the function class $\Theta$
which does not determine $\theta$ up to finitely many real-valued
parameters. Thus we consider a nonparametric estimation problem.
Moreover we assume that $E X_1 = 0$ and $P[\|X_1\|_2 \geq x]   \leq
C_{X,0} \exp(-C_{X,1} x^{C_{X,2}})$ for all $x>0$ and some finite
constants $C_{X,0},C_{X,1},C_{X,2}>0$ where $\|\cdot\|_p$, $p\geq1$
denotes the $L_p([0,1])$-norm of some element of that space. Thus the
tails of the design distribution are restricted. Such conditions are
usual in nonparametric regression problems.

The FLR model has obtained considerable attention in the statistical
community during the last years, which is reflected in the large amount
of literature on this topic. Various of estimation procedures have been
proposed to make the regression function $\theta$ empirically
accessible (see, e.g.,~\cite{Cardotetal1999,Cardotetal2003,CardotandJohannes2009,Crambesetal2009,Delaigleetal2009}). The minimax convergence rates in FLR are
investigated, for example, in~\cite{CaiandHall2006,HallandHorowitz2007,CardotandJohannes2009}. In \cite
{CaiandZhou2008}, adaptive estimation in FLR is considered.
Generalizations of FLR are discussed in \cite
{MuellerandStadtmueller2005}. A central limit theorem for FLR is
derived in~\cite{Cardotetal2007}. In~\cite{Yaoetal2005}, practical
applications of FLR in the field of medical statistics are described;
the authors consider two real data sets on primary biliary cirrhosis
and systolic blood pressure. For a comprehensive introduction to the
field of functional data analysis in general, see~\cite{RamsayandSilverman2005}.

In order to compare two statistical models, it is useful to prove
asymptotic equivalence between those models. For the basic concept and
a detailed description of this strong asymptotic property, we refer to
\cite{LeCam1964} and~\cite{LeCamandYang2000}. Also, a review on this
topic is given in the following section. As an important feature, if
two models $\mathfrak{E}_{1,n}$ and $\mathfrak{E}_{2,n}$ are
asymptotically equivalent, then $\mathfrak{E}_{1,n}$ adopts optimal
convergence rates and sharp asymptotic constants with respect to \textit{any} bounded loss function from model $\mathfrak{E}_{2,n}$ and vice
versa. Thus, the theory of asymptotic equivalence does not only capture
special loss functions such as the mean integrated squared error (MISE)
or the pointwise mean squared error (MSE) but includes various types of
semi-metrics between the estimator and the target function $\theta$ and
also addresses the estimation of characteristics of $\theta$, such as
its support or its mode. Furthermore, superefficiency phenomena also
coincide in both models when considering subclasses $\Theta'$ of the
target parameter space $\Theta$. In particular, research has focussed
on proofs of asymptotic equivalence of experiments where $n$ i.i.d.
data are observed, whose distribution depends on some parameter $\theta
\in\Theta$, and experiments where $\theta$ occurs in the drift of an
empirically accessible It\^{o} process. For instance, Nussbaum \cite
{Nussbaum1996} considers an asymptotically equivalent white noise model
for density estimation, while Brown and Low~\cite{BrownandLow1996}
introduce such a model for nonparametric regression. In recent related
literature on regression problems, Carter~\cite{Carter2007} studies the
case of unknown error variance, and Reiss~\cite{Reiss2008} extends
asymptotic equivalence to the multivariate setting.

Returning to model (\ref{eq:1.1}), we suppose that the nuisance
parameters $\sigma$ and $P_X$, that is, the distribution of the $X_j$,
are known. That allows us to exclude those quantities from the
parameter space of the experiment and to fully concentrate on the
estimation of $\theta$. This condition is also imposed in most papers
dealing with asymptotic equivalence for nonparametric regression
experiments. The work of~\cite{Carter2007} represents an exception
where the corresponding white noise model becomes more difficult and,
apparently, less useful to derive adoptable asymptotic properties. With
respect to asymptotic equivalence, we restrict our consideration to the
case of known $P_X$. However, in Section~\ref{S:new}, we will show that
the sharp minimax asymptotics with respect to the MISE are extendable
to the case of unknown design distribution.

The main purpose of the current work is to prove asymptotic equivalence
of model (\ref{eq:1.1}) and a statistical inverse problem in the white
noise setting. That latter model is described by the observation of an
It\^{o} process $Y(t)$, $t\in[0,1]$, $Y(0)=0$, driven by the
stochastic differential equation
%
%e1.2 ###
\begin{equation} \label{eq:1.1.1}
dY(t)   =    [K \theta ](t)\,dt + n^{-1/2} \sigma \,dW(t) ,
\end{equation}
 where $W(t)$ denotes a standard Wiener process on the
interval $[0,1]$, and $K$ denotes a linear operator mapping from the
Hilbert space $L_2([0,1])$ to itself. These models are also widely
studied in mathematical statistics (see, e.g.,~\cite
{CavalierandTsybakov2002} and~\cite{Goldenshlugeretal2006}). They have
their applications in the field of signal deblurring and econometrics.
We will concentrate on a specific version of model (\ref{eq:1.1.1})
where $K$ is equal to the unique positive symmetric square root $\Gamma
^{1/2}$ of the covariance operator $\Gamma$, that is, $\Gamma^{1/2}
\Gamma^{1/2} = \Gamma$ and $\Gamma f   =   \int E X_1(\cdot) X_1(t)
f(t)\,dt$ for any $f \in L_2([0,1])$. Thus, the observation $Y(t)$,
$t\in[0,1]$, is defined by $Y(0)=0$ and
%
%e1.3 ###
\begin{equation} \label{eq:1.1.2}
dY(t)   =    [\Gamma^{1/2} \theta ](t)\,dt + n^{-1/2} \sigma
\,dW(t) .
\end{equation}
 In~\cite{CardotandJohannes2009}, the authors remark on the
similarity of models (\ref{eq:1.1}) and (\ref{eq:1.1.2}). In the
current paper, we will rigorously establish asymptotic equivalence
between those models. As an interesting feature, additional observation
of the data $X_1,\ldots ,X_n$ would be redundant in model (\ref
{eq:1.1.2}). All information about the design points is recorded by
$\Gamma$ in (\ref{eq:1.1.2}). Therefore, all what is observed in the
corresponding white noise experiment is the process $Y(t)$, $t\in
[0,1]$. After the general introduction to the property of asymptotic
equivalence as used in the current paper in Section~\ref{S1.0}, we will
first prove (nonasymptotic) equivalence of model (\ref{eq:1.1}) and an
empirical  version of model (\ref{eq:1.1.2}) where $\Gamma$ is replaced by
a noisy counterpart in Section~\ref{S2}. In Section~\ref{S3}, we prove
asymptotic equivalence of (\ref{eq:1.1}) and (\ref{eq:1.1.2}) under
some additional technical conditions. In Section~\ref{S:new}, we show
that the sharp lower bound which follows from the results of the
previous section can be attained by specific estimators in the
realistic case of unknown design distribution. A discussion of the
findings and their conclusions are provided in Section~\ref{S5}.

%s2 ###
\section{Asymptotic equivalence} \label{S1.0}

To recall the definition of asymptotic equivalence, we consider two
(sequences of) statistical experiments $\mathfrak{E}_{j,n} = (\Omega
_{j,n},\break\mathfrak{A}_{j,n}, P_{j,n,\theta})$, $j=j_1,j_2$, with a joint
parameter space $\Theta$, which may depend on $n$. The LeCam distance
between $\mathfrak{E}_{j_1,n}$ and $\mathfrak{E}_{j_2,n}$ is defined by
\[
\Delta(\mathfrak{E}_{j_1,n},\mathfrak{E}_{j_2,n})   =   \max_{k=1,2}
\inf_{K\in\mathfrak{K}_{j_k,n}} \sup_{\theta\in\Theta} \|
K(P_{j_k,n,\theta}) - P_{j_{3-k},n,\theta}\|_{\mathrm{TV}} ,
\]
 where $\|\cdot\|_{\mathrm{TV}}$ is the total variation distance, and
$\mathfrak{K}_{j_k,n}$ denotes the collection of so-called transitions
(see~\cite{Strasser1985} and~\cite{Nussbaum1996} for their exact
definition). The statistical experiments $\mathfrak{E}_{j_1,n}$ and
$\mathfrak{E}_{j_2,n}$ are called asymptotically equivalent if $\Delta
(\mathfrak{E}_{j_1,n},\mathfrak{E}_{j_2,n})$ converges to zero as $n\to
\infty$, while they are called equivalent if $\Delta(\mathfrak
{E}_{j_1,n},\mathfrak{E}_{j_2,n})=0$ for all $n$.

In the framwork of our note, we will not use that general definition of
(asymptotic) equivalence but our proofs lean on following sufficient
conditions for these properties:
\begin{longlist}[(iii)]
 \item[(i)] We consider the following sufficient condition for
asymptotic equivalence of $\mathfrak{E}_{j_1,n}$ and $\mathfrak
{E}_{j_2,n}$: We define the sets ${\cal R}_{j,n,\theta}$, $j=j_1,j_2$,
$\theta\in\Theta$, which contains all real-valued integrable random
variables $R$ on the domain $\Omega_{j,n}$ satisfying $|R|\leq1$ a.s.
Thus any kind of bounded loss functions are captured by the classes
${\cal R}_{j,n,\theta}$ so that the expectation $E R$ with respect to
the distribution $P_{j,n,\theta}$ describes an arbitrary bounded and
normalized statistical risk for estimating the parameter $\theta$ under
the observation scheme $\mathfrak{E}_{j,n}$. Now we define two
sequences $(T_{j_k,j_{3-k},n})_n$, $k=1,2$, of $(\mathfrak
{A}_{j_k},\mathfrak{A}_{j_{3-k}})$-measurable mappings from $\Omega
_{j_k}$ to $\Omega_{j_{3-k}}$. As an essential condition, the
$T_{j_k,j_{3-k},n}$ must not depend on $\theta$. Hence,
$T_{j_k,j_{3-k},n}$ may be interpreted as a transformation of the data
from an observation contained in the space $\Omega_{j,n}$ to an
observation which lies in $\Omega_{3-j,n}$. Thus a statistician who
intends to construct an estimation procedure for $\theta$ may always
apply this transformation $T_{j_k,j_{3-k},n}$ to an observation $\omega
\in\Omega_{j_k,n}$. Then we obtain asymptotic equivalence of $\mathfrak
{E}_{j_1,n}$ and $\mathfrak{E}_{j_2,n}$ when we can show the existence
of such transformation sequences $(T_{j_k,j_{3-k},n})_n$, $k=1,2$, so that
%
%e2.1 ###
\begin{equation} \label{eq:aseq}\qquad
\sup_{\theta\in\Theta}  \sup_{R_{j_{3-k},n,\theta} \in{\cal
R}_{j_{3-k},n,\theta}}   | E R_{j_{3-k},n,\theta}   -   E
(R_{j_{3-k},n,\theta} \circ T_{j_k,j_{3-k},n})  |   \longrightarrow
  0
\end{equation}
 as $n\to\infty$ for all $k=1,2$. Accordingly, we have
equivalence if the left-hand side in (\ref{eq:aseq}) equals $0$ for any
$n$. Intuitively speaking, after transforming the data drawn from model
$\mathfrak{E}_{j_1,n}$ according to $T_{j_1,j_2,n}$, the distance
between any bounded statistical risk in model $\mathfrak{E}_{j_2,n}$ on
one hand and for the transformed data from model $\mathfrak{E}_{j_1,n}$
becomes small for large $n$ or is equal to zero for any~$n$,
respectively. The same condition must also hold true when exchanging
$j_1$ and $j_2$.

In the specific framework of our note, we assume, in addition, that all
transformations $T_{j_k,j_{3-k},n}$ must not depend on the nuisance
parameter $\sigma$. That compensates the unrealistic condition of known
$\sigma$. In particular, $\sigma$ is not used to transform the data or
to construct decision procedures or estimators. Therefore, our results
also addresses the case of unknown $\sigma$. Nevertheless, $\sigma$
must be viewed as uninteresting for the statistician, that is, it must
not explicitly occur in the loss functions $R_{j_k,j_{k+1}}$. Thus, the
problem of estimating $\sigma$ is not covered by our approach.

 \item[(ii)]  Assume that the experiment $\mathfrak{E}_{j_2,n}$
describes the observation of $T(\omega)$ for $\omega\in\Omega
_{j_1,n,\theta}$ in experiment $\mathfrak{E}_{j_1,n}$ where $T$ is a
sufficient statistic for $\theta$ in experiment $\mathfrak{E}_{j_1,n}$.
Then $\mathfrak{E}_{j_1,n}$ and $\mathfrak{E}_{j_2,n}$ are
statistically equivalent (i.e., their LeCam distance vanishes) and,
hence, asymptotically equivalent. That assertion holds true whenever
the experiments are Polish spaces. This criterion is satisfied as all
probability spaces considered in the current work are $\mathbb{R}^d$,
$C([0,1])$, $L_2([0,1])$ and some set products of those classes (see
Lemma 3.2 in~\cite{BrownandLow1996}).

 \item[(iii)]  If some experiments $\mathfrak{E}_{j_1,n}$ and
$\mathfrak{E}_{j_2,n}$ on one hand, and $\mathfrak{E}_{j_2,n}$ and
$\mathfrak{E}_{j_3,n}$ on the other hand, are (asymptotically)
equivalent, then $\mathfrak{E}_{j_1,n}$ and $\mathfrak{E}_{j_3,n}$ are
(asymptotically) equivalent, too. Also, (asymptotic) equivalence of
$\mathfrak{E}_{j_1,n}$ and $\mathfrak{E}_{j_2,n}$ is a symmetric
relation between the experiments.

 \item[(iv)]  Assume that some experiments $\mathfrak{E}_{j_1,n}$
and $\mathfrak{E}_{j_2,n}$ may be decomposed into two independent
experiments $\mathfrak{E}_{j_1,n,k}$ and $\mathfrak{E}_{j_2,n,k}$,
$k=1,2$, respectively. Moreover, we suppose that the experiments
$\mathfrak{E}_{j_1,n,1}$ and $\mathfrak{E}_{j_2,n,1}$ on one hand and
the experiments $\mathfrak{E}_{j_1,n,2}$ and $\mathfrak{E}_{j_2,n,2}$
on the other hand are (asymptotically) equivalent. Then, the combined
experiments $\mathfrak{E}_{j_1,n}$ and $\mathfrak{E}_{j_2,n}$ are
(asymptotically) equivalent as well.
\end{longlist}

Now, $\mathfrak{E}_{1,n}$ denotes the underlying experiment of the FLR
model (\ref{eq:1.1}); it is defined by $\Omega_{1,n} = C([0,1])^{(n)}
\times\mathbb{R}^{(n)}$, $\mathfrak{A}_{1,n}$ denotes the Borel $\sigma
$-algebra when considering the uniform metric on the functional
components and the Euclidean \mbox{metric} on the real-valued components of
$\Omega_{1,n}$. The corresponding probability measures $P_{1,n,\theta}$
are well defined by the assumptions of the model~(\ref{eq:1.1}). The
parameter space $\Theta\subseteq L_2([0,1])$ will be specified later.
Still, the observations $(\mathbf{X},\mathbf{Y})$ may be viewed as random
variables having their domain on some basic probability space $(\Omega
,\mathfrak{A},P)$.

%s3 ###
\section{Empirical covariance operator} \label{S2}

We define the linear covariance operator $\Gamma\dvtx  L_2([0,1]) \to
L_2([0,1])$ by
\[
\Gamma f   =   \int E X_1(\cdot) X_1(t) f(t)\,dt   \qquad\forall f
\in L_2([0,1]) .
\]
 Writing $K(s,t) = E X_1(s) X_1(t)$, we realize that $\Gamma$
is a Hilbert--Schmidt integral operator where
\[
\int_0^1 \int_0^1 |K(s,t)|^2\,ds \,dt   \leq  (E \|X_1\|_2^2)^2   <
\infty ,
\]
 by the Cauchy--Schwarz inequality and the tail condition
imposed on the distribution of $\|X_1\|_2$. Hence $\Gamma$ is a
continuous and compact operator. We have $K(s,t) = K(t,s)$ for all $s,t
\in[0,1]$ so that the operator $\Gamma$ is self-adjoint. Furthermore,
it is positive; that is, by Fubini's theorem we have
\[
\langle f, \Gamma f\rangle  =   E | \langle X_1 , f \rangle
|^2   \geq  0
\]
 for any $f\in L_2([0,1])$.\vadjust{\goodbreak}

Then, some well-known results from functional analysis, in particular
spectral theory for compact operators, may be applied. There exists an
orthonormal basis $\{\varphi_j\}_{j\geq1}$ of the separable Hilbert
space $L_2([0,1])$ which consists of eigenfunctions of $\Gamma$. The
corresponding eigenvalues are denoted by $\lambda_j\geq0$. The
sequence $(\lambda_n)_n$ converges to zero and may be viewed as
monotonously decreasing without loss of generality. Those results are
also used, for example, in~\cite{CaiandHall2006}. Furthermore, for
$\Gamma$ as for any compact self-adjoint positive operator from
$L_2([0,1])$ to itself, there exists a unique compact self-adjoint
positive operator $\Gamma^{1/2}$ from $L_2([0,1])$\vspace*{1pt} to itself such that
$(\Gamma^{1/2})^2=\Gamma$; then $\Gamma^{1/2}$ is called the square
root of $\Gamma$. We have $\Gamma^{1/2} \varphi_j   =   \lambda
_j^{1/2} \varphi_j$ for any $j\geq1$.

We may define an empirically accessible version $\hat{\Gamma}$ of
$\Gamma$ by replacing the expectation by the average; more precisely,
we have
\[
\hat{\Gamma} f   =   \frac1n \sum_{j=1}^n \int X_j(\cdot) X_j(t) f(t)\,dt   \qquad\forall f \in L_2([0,1]) .
\]
 Thus, $\hat{\Gamma}$ may be viewed as the operator $\Gamma$
when $P_X$ equals the uniform distribution on the discrete set $\{
X_1,\ldots ,X_n\}$. Therefore, all properties derived for $\Gamma$ in
the previous paragraph can be taken over to $\hat{\Gamma}$. In
particular, $\hat{\varphi}_j$, integer $j\geq1$, denotes the
orthonormal basis of the eigenfunctions of $\hat{\Gamma}$ with the
eigenvalues $\hat{\lambda}_j$.

Now we consider the conditional probability density $p_{1,n,\theta
}(y_1,\ldots , y_n \mid X_1,\ldots ,X_n)$ of the data $Y_1,\ldots ,Y_n$
given the design functional observations $X_1,\ldots ,X_n$ in model (\ref
{eq:1.1}). This density shall be understood with respect to the
$n$-dimensional Lebesgue measure. We derive that
%
%e3.1 ###
\begin{eqnarray}\label{eq:S2.n.1}
&&p_{1,n,\theta}(y_1,\ldots ,y_n \mid X_1,\ldots ,X_n) \nonumber\\
&& \qquad    =
  (2\pi)^{-n/2} \sigma^{-n} \prod_{j=1}^n \exp \biggl( - \frac
1{2\sigma^2} (y_j - \langle X_j,\theta\rangle)^2 \biggr) \\
&& \qquad    =   (2\pi)^{-n/2} \sigma^{-n}   \exp \bigl( - \|\mathbf{y} - \mathbf{x}\|
^2 / (2\sigma^2) \bigr) , \nonumber
\end{eqnarray}
 with the vectors $\mathbf{y} = (y_1,\ldots ,y_n)^T$ and $\mathbf{x}
= (\langle X_1,\theta\rangle,\ldots , \langle X_n,\theta\rangle)^T$.
Moreover, \mbox{$\|\cdot\|$} denotes the Euclidean norm. Expanding $\theta\in
\Theta\subseteq L_2([0,1])$ in the orthonormal basis $\{\hat{\varphi
}_j\}_{j\geq1}$ gives us that
%
%e3.2 ###
\begin{equation} \label{eq:S2.n.2} \langle X_j , \theta\rangle  =
\sum_{k=1}^\infty\langle X_j , \hat{\varphi}_k \rangle  \langle\hat
{\varphi}_k , \theta\rangle .
\end{equation}
 We impose the following condition on the distribution $P_X$:
%
%e3.3 ###
\begin{eqnarray}\label{eq:Cond_X}
 \qquad P[X_1 \in L] &=& 0 ,   \qquad  \mbox{for any deterministic linear subspace } \nonumber
\\[-8pt]
\\[-8pt]
  &&     \hspace*{32pt}L
\subseteq L_2([0,1])\mbox{ with }\dim L <\infty .
\nonumber
\end{eqnarray}
 Intuitively, this assumption provides that the probability
mass of the $X_j$ fills the whole of $L_2([0,1])$. Somehow, (\ref
{eq:Cond_X}) is the functional data analog for continuity of a
distribution of some real-valued random variables. It is satisfied when
we take an appropriate Gaussian process for $X_1$, for instance.
Condition (\ref{eq:Cond_X}) yields that the linear space generated by
$X_1,\ldots ,X_n$ is $n$-dimensional almost surely. Otherwise, at least
one of the $X_j$ must be included in the linear hull of the other
design variables. According to (\ref{eq:Cond_X}) that occurs with
probability zero when employing the conditional probability measure
given the data $X_1,\ldots ,X_{j-1},X_{j+1},\ldots ,X_n$. Finally,
applying the expectation, we obtain the desired result for the
unconditional distribution.

We realize that the range of $\hat{\Gamma}$ is included in the linear
hull of $X_1,\ldots ,X_n$. By definition, $\hat{\varphi}_j$ is contained
in that $n$-dimensional space whenever $\hat{\lambda}_j > 0$. As the
$\hat{\varphi}_j$ form an orthonormal basis at most $n$ of the
eigenvalues $\hat{\lambda}_j$ are nonvanishing. Furthermore, the
linear independence of the $X_1,\ldots ,X_n$ implies that the functions
$\hat{\Gamma} X_k   =   n^{-1} \sum_{j=1}^n \langle X_j,X_k \rangle
  X_j$, $k=1,\ldots ,n$, are linearly independent, too, so that the
range of $\hat{\Gamma}$ is equal to the linear hull of $X_1,\ldots
,X_n$. Clearly, the range of $\hat{\Gamma}$ also coincides with the
linear hull of all $\hat{\varphi}_j$ with $\hat{\lambda}_j > 0$, from
what follows $\langle X_j , \hat{\varphi}_k \rangle= 0$ for all
$j=1,\ldots ,n$ and $k>n$. Also, we have $\hat{\lambda}_j > 0$ for
$j=1,\ldots ,n$ and $\hat{\lambda}_j = 0$ for $j>n$. Hence, (\ref
{eq:S2.n.2}) leads to the representation
%
%e3.4 ###
\begin{equation} \label{eq:S2.n.3}
\langle X_j , \theta\rangle  =   \sum_{k=1}^n \langle X_j , \hat
{\varphi}_k \rangle  \langle\hat{\varphi}_k , \theta\rangle
\end{equation}
 for all $j=1,\ldots ,n$. Equation (\ref{eq:S2.n.3}) is
equivalent to the system of linear equations $\mathbf{x} = \mathbf{Q} \mathbf{f}$ with the vector $\mathbf{f} = (\langle\hat{\varphi}_1 , \theta\rangle
,\ldots , \langle\hat{\varphi}_n , \theta\rangle)^T$ and the matrix
$\mathbf{Q}$ with the components $Q_{j,k} = \langle X_j , \hat{\varphi}_k
\rangle$, $j,k = 1,\ldots ,n$. Then the conditional density
$p_{1,n,\theta}$ as in (\ref{eq:S2.n.1}) may be written as
%
%e3.5 ###
\begin{equation} \label{eq:S2.n.4}
  p_{1,n,\theta}(y_1,\ldots ,y_n \mid X_1,\ldots ,X_n)   =   (2\pi
)^{-n/2} \sigma^{-n}   \exp\bigl ( - \|\mathbf{y} - \mathbf{Q} \mathbf{f}\|^2 /
(2\sigma^2) \bigr) .\hspace*{-30pt}
\end{equation}
 We consider that the $(k,k')$th component of the matrix $\mathbf{Q}^T \mathbf{Q}$ is equal to
\[
\sum_{j=1}^n \langle X_j,\hat{\varphi}_k \rangle  \langle X_j,\hat
{\varphi}_{k'} \rangle  =   n \langle\hat{\Gamma} \hat{\varphi}_k ,
\hat{\varphi}_{k'}\rangle  =   n \hat{\lambda}_k \cdot\delta
_{k,k'} .
\]
 Thus $\mathbf{Q}^T \mathbf{Q}$ is a diagonal matrix containing $n
\hat{\lambda}_{k}$ as its $(k,k)$th component. We denote the diagonal
matrix having $n^{1/2} \hat{\lambda}_{k}^{1/2}$ as its $(k,k)$th
component by $\mathbf{D}$. Obviously, $\mathbf{D}$ is invertible, and we
define $\mathbf{A} = \mathbf{Q} \mathbf{D}^{-1}$. We have
\[
\mathbf{A}^T \mathbf{A}   =   \mathbf{D}^{-1} \mathbf{Q}^T \mathbf{Q} \mathbf{D}^{-1}
  =   \mathbf{I} ,
\]
 where $\mathbf{I}$ denotes the identity matrix. Also, this
yields that $\mathbf{A} \mathbf{A}^T = \mathbf{I}$ and that $\mathbf{A}$ is an
orthogonal matrix. Thus, $\|\mathbf{A} \mathbf{v}\|   =   \|\mathbf{v}\|$ for
any vector $\mathbf{v} \in\mathbb{R}^n$. Equality (\ref{eq:S2.n.4})
provides that
%
%e3.6 ###
\begin{eqnarray} \label{eq:S2.n.5}
&&p_{1,n,\theta} (y_1,\ldots ,y_n \mid X_1,\ldots ,X_n)\nonumber \\
 && \qquad
=   (2\pi)^{-n/2} \sigma^{-n}   \exp \bigl( - \|\mathbf{A} \mathbf{A}^T \mathbf{y} - \mathbf{A} \mathbf{D} \mathbf{f}\|^2 / (2\sigma^2) \bigr) \\
&& \qquad    =   (2\pi)^{-n/2} \sigma^{-n}   \exp \bigl( - \|\mathbf{A}^T \mathbf{y}
- \mathbf{D} \mathbf{f}\|^2 / (2\sigma^2) \bigr) .\nonumber
\end{eqnarray}

Referring to the notation of (\ref{eq:aseq}), we consider the
expectation $E R_{1,n,\theta}(\mathbf{X},\break\mathbf{Y})$ where $R_{1,n,\theta}
\in{\cal R}_{1,n,\theta}$. We derive that
%
%e3.7 ###
\begin{eqnarray}
\label{eq:S2.n.6}
  E R_{1,n,\theta}(\mathbf{X},\mathbf{Y})   &=&   E   E  (R_{1,n,\theta
}(\mathbf{X},\mathbf{Y}) \mid\mathbf{X} )\nonumber \\
    &= &  E \idotsint R_{1,n,\theta}(X_1,\ldots ,X_n;y_1,\ldots
    ,y_n)\nonumber\\
&&\hphantom{E \idotsint}{}\times p_{1,n,\theta}(y_1,\ldots ,y_n \mid X_1,\ldots ,X_n) \,  dy_1\cdots dy_n \nonumber
\\[-8pt]
\\[-8pt]
   &=&   E \idotsint R_{1,n,\theta}(X_1,\ldots ,X_n;\mathbf{A}\mathbf{z})
(2\pi)^{-n/2} \nonumber\\
&&\hphantom{E \idotsint}
{}\times\sigma^{-n}   \exp \bigl( - \|\mathbf{z} - \mathbf{D} \mathbf{f}\|
^2 / (2\sigma^2) \bigr) \,  dz_1 \cdots dz_n \nonumber\\
   &=&   E R_{1,n,\theta}(\mathbf{X},\mathbf{A}\mathbf{Z}) ,
\nonumber
\end{eqnarray}
 where $\mathbf{Z} = (Z_1,\ldots ,Z_n)^T$ denotes a vector
consisting of independent normally distributed random variables where
$Z_k$ has the mean
\[
n^{1/2} \hat{\lambda}_k^{1/2} \langle\hat{\varphi}_k , \theta\rangle
  =   \langle\hat{\varphi}_k , n^{1/2} \hat{\Gamma}^{1/2} \theta
\rangle ,
\]
 and the variance $\sigma^2$, conditionally on the $\sigma
$-algebra generated by $\mathbf{X}$. Therefore, the $Z_k$ may be
represented as
%
%e3.8 ###
\begin{equation} \label{eq:nn.1}
Z_k   =   \langle\hat{\varphi}_k ,
n^{1/2} \hat{\Gamma}^{1/2} \theta\rangle  +   \sigma\varepsilon_k
, \qquad k=1,\ldots ,n ,
\end{equation}
 where $\varepsilon_1,\ldots ,\varepsilon_n$ are i.i.d.
$N(0,1)$-distributed random variables. The $\varepsilon_j$ are
independent of the $\sigma$-algebra generated by $\mathbf{X}$. We have
applied the integral transformation $\mathbf{y} = \mathbf{A} \mathbf{z}$ where
$\det\mathbf{A} = \pm1$ due to the orthogonality of $\mathbf{A}$. Note that
the sign of the eigenfunctions $\hat{\psi}_j$ may still be chosen; we
can arrange that $\det\mathbf{A} = 1$.

Now we define the statistical experiment $\mathfrak{E}_{2,n}$ with the
same parameter space $\Theta$ as $\mathfrak{E}_{1,n}$, $(\Omega
_{2,n},\mathfrak{A}_{2,n}) = (\Omega_{1,n},\mathfrak{A}_{1,n})$ and
$P_{2,n,\theta}$ as the probability measure generated by the random
variable $(\mathbf{X},\mathbf{Z})$ with $\mathbf{Z}$ as in (\ref{eq:S2.n.6}). In
the notation of Section~\ref{S1.0}, paragraph (i), we use the mapping
$T_{2,1,n}\dvtx  \Omega_{2,n} \to\Omega_{1,n}$ defined by $T_{2,1,n}(\mathbf{x},\mathbf{z}) = (\mathbf{x},\mathbf{A}\mathbf{z})$, $\mathbf{x} \in
C_0([0,1])^{(n)}$, $\mathbf{z} \in\mathbb{R}^n$, as the data
transformation from $\mathfrak{E}_{2,n}$ to $\mathfrak{E}_{1,n}$. By
definition, the matrix $\mathbf{A}$ does not depend on the parameter
$\theta$ but only on the data $X_1,\ldots ,X_n$ and the known
orthonormal basis $\{\hat{\varphi}_j\}_{j\geq1}$. Also, it does not
depend on $\sigma$ as requested in the previous section. We have
already derived that $\mathbf{A}$ is an orthogonal matrix so that
$T_{2,1,n}$ is a bijective\vadjust{\goodbreak} mapping from the set $C_0([0,1])^{(n)}
\times\mathbb{R}^n$ to itself. Hence, its reverse mapping
$T_{2,1,n}^{-1}$ may be used as the data transformation $T_{1,2,n}$.
Then, according to (\ref{eq:aseq}), we have proved the following lemma.

\begin{lem} \label{L:1}
Under condition (\ref{eq:Cond_X}), the statistical experiments
$\mathfrak{E}_{1,n}$ and $\mathfrak{E}_{2,n}$ are equivalent.
\end{lem}

The random variables $\varepsilon_j$, integer $j$, as occurring in (\ref
{eq:nn.1}), may be represented by
\[
\varepsilon_j   =   \int_0^1 \hat{\varphi}_j(t) \, dW(t) ,
\]
 where $W$ denotes a standard Wiener process on $[0,1]$ which
is independent of $\mathbf{X}$. We deduce that the $\varepsilon
_1,\varepsilon_2,\ldots $ are an independent sequence of
$N(0,1)$-distributed random variables. Moreover, they are independent
of $X_1,\ldots ,X_n$ although $\hat{\varphi}_j$ depends on these design
variables. That can be shown via the conditional characteristic
function of $(\varepsilon_1,\varepsilon_2,\ldots )$ given $X_1,\ldots
,X_n$; that is,
\begin{eqnarray*}
E  \Biggl[\exp \Biggl( i \sum_{j=1}^\infty\int s_j \hat
{\varphi}_j(t)\,dW(t) \Biggr) \Bigm| X_1,\ldots ,X_n \Biggr]
 &   =&   \exp \Biggl(-
\frac12 \sigma^2  \Biggl\|\sum_{j=1}^\infty s_j \hat{\varphi}_j  \Biggr\|
_2^2 \Biggr) \\
 &   =&   \exp \Biggl(- \frac12 \sigma^2 \sum_{j=1}^\infty
s_j^2 \Biggr)
\end{eqnarray*}
 for all real-valued sequences $(s_m)_{m\geq1}$ with $\sum
_{m=1}^\infty s_m^2 < \infty$. Applying the expectation to the above
equality, the unconditional characteristic function of $(\varepsilon
_1,\varepsilon_2,\ldots )$ turns out to coincide with the conditional
one. We have
%
%e3.9 ###
\begin{equation} \label{eq:S2.n.6.1}
Z_j   =   \langle\hat{\varphi
}_j , n^{1/2} \hat{\Gamma}^{1/2} \theta\rangle  +   \sigma\int_0^1
\hat{\varphi}_j(t)\,dW(t)   =   \int \hat{\varphi}_j \, dZ(t)
\end{equation}
 for all $j=1,\ldots ,n$ where $Z(t)$, $t\in[0,1]$, denotes an
It\^{o} process satisfying
%
%e3.10 ###
\begin{equation} \label{eq:S2.n.7}
dZ(t)   =   n^{1/2}  [\hat{\Gamma}^{1/2} \theta ](t)\,dt  +
\sigma\, dW(t) ,
\end{equation}
 and $Z(0)=0$. The differential $dZ(t)$ shall be understood in
the It\^{o} sense.

Now we define the statistical experiment $\mathfrak{E}_{3,n}$ with a
completely functional observation structure. We fix that $\Omega_{3,n}
= C([0,1])^{(n+1)}$ with the corresponding Borel $\sigma$-algebra
$\mathfrak{A}_{3,n}$. The probability measure $P_{3,n,\theta}$ is
defined via the observation of $\mathbf{X}$ as in $\mathfrak{E}_{2,n}$ and
the It\^{o} process $Z(t)$, $t\in[0,1]$, as defined in~(\ref
{eq:S2.n.7}). The definition (\ref{eq:S2.n.6.1}) of $Z_j$ can be
extended to $j>n$ straightforwardly. As $\hat{\lambda}_j=0$, we obtain that
\[
Z_j    =   \sigma\int_0^1 \hat{\varphi}_j(t)\,dW(t)   \qquad\forall
j>n .
\]
 Moreover, $Z(t)$ is uniquely determined by the $Z_j$ for all
integers $j\geq1$ and vice versa. That can be seen as follows:
\[
Z(t) = \int1_{[0,t]}(s)\,dZ(s)   =   \sum_{j=1}^\infty\langle
1_{[0,t]} , \hat{\varphi}_j \rangle  Z_j
\]
 for all $t\in[0,1]$ where the infinite sum must be
understood as an $E\|\cdot\|_2^2$-limit. That seems to cause some
troubles as we only observe one element of the probability space.
However, convergence in probability implies almost sure convergence of
a subsequence so that $Z(t)$ is fully accessible by the observation of
all $Z_j$. On the other hand, by a similar argument, all $Z_j$ are
accessible (in practice, that means approximable arbitrarily precisely)
by a trajectory of the process $Z$.

Hence the data set $\{Z_j   \dvtx    j > n\}$ is independent of the
$Z_1,\ldots ,Z_n$, conditionally on the $\sigma$-algebra generated by
$\mathbf{X}$. Furthermore, the distribution of the $Z_j$, $j > n$, does
not depend on the target parameter $\theta$ so that $Z_j$, for $j>n$,
does not contain any information about $\theta$. We conclude that
$(\mathbf{X},Z_1,\ldots ,Z_n)$ is a sufficient statistic for the
observation scheme in the experiment $\mathfrak{E}_{3,n}$. We can
utilize result (ii) from Section~\ref{S1.0} in order to prove
equivalence of the experiments $\mathfrak{E}_{2,n}$ and $\mathfrak
{E}_{3,n}$. Considering paragraph (iii) from Section~\ref{S1.0}, we may
establish equivalence of the experiments $\mathfrak{E}_{1,n}$ and
$\mathfrak{E}_{3,n}$. This result is presented in the following theorem.

\begin{teo} \label{T:1}
Under condition (\ref{eq:Cond_X}), the FLR statistical experiment
$\mathfrak{E}_{1,n}$ is equivalent to the model $\mathfrak{E}_{3,n}$
where one observes $\mathbf{X}$ and the It\^{o} process $Z(t)$, $t \in
[0,1]$, as defined in (\ref{eq:S2.n.7}).
\end{teo}

%s4 ###
\section{Asymptotic approximation} \label{S3}

In the previous section, we have derived a statistically equivalent
white noise model for the FLR problem. However, the It\^{o} process $Z$
in (\ref{eq:S2.n.6.1}) contains the noisy operator $\hat{\Gamma}$ in
its construction. In the current section, we will replace it by the
covariance operator $\Gamma$.

For that purpose, we split the original experiment $\mathfrak{E}_{1,n}$
into two independent parts $\mathfrak{E}_{1,n,1}$ and $\mathfrak
{E}_{1,n,2}$ where $\mathfrak{E}_{1,n,1}$ is based on the observation
of the data $(X_j,Y_j)$, $j=1,\ldots ,m$, and $\mathfrak{E}_{1,n,2}$
consists of the residual data $(X_j,Y_j)$, $j=m+1,\ldots ,n$. The
selection of the integer parameter $m$ is deferred. The strategy of
splitting the sample in the current context leans on \cite
{Nussbaum1996}. Applying Theorem~\ref{T:1} to each of the experiments
$\mathfrak{E}_{1,n,k}$, $k=1,2$, we obtain equivalence of $\mathfrak
{E}_{1,n,k}$ and the experiments $\mathfrak{E}_{4,n,k}$ for $k=1,2$.
Therein, $\mathfrak{E}_{4,n,1}$ is defined by the observation of $\mathbf{X}_1 = (X_1,\ldots ,X_m)$ and the It\^{o} process $Z_1(t)$, $t \in
[0,1]$, specified by $Z_1(0)=0$ and
\[
dZ_1(t)   =   m^{1/2}  [\hat{\Gamma}_1^{1/2} \theta ](t)\,dt
+   \sigma\, dW_1(t) ,
\]
 and accordingly $\mathfrak{E}_{4,n,2}$ is defined by the
observation of $\mathbf{X}_2 = (X_{m+1},\ldots ,X_n)$ and the It\^{o}
process $Z_2(t)$, $t \in[0,1]$, specified by $Z_2(0)=0$ and
\[
dZ_2(t)   =   (n-m)^{1/2}  [\hat{\Gamma}_2^{1/2} \theta ](t)\,dt
  +   \sigma\, dW_2(t) .
\]
 Furthermore, $\hat{\Gamma}_k$, $k=1,2$, denotes the empirical
covariance operator constructed by the data $\mathbf{X}_1$ and $\mathbf{X}_2$, respectively. Also note that $W_1$ and $W_2$ are two independent
standard Wiener processes. Using criterion (iv) in Section~\ref{S1.0},
the experiment $\mathfrak{E}_{4,n}$, which combines the independent
experiments $\mathfrak{E}_{4,n,1}$ and $\mathfrak{E}_{4,n,2}$, we
deduce that $\mathfrak{E}_{4,n}$ and $\mathfrak{E}_{1,n}$ are equivalent.

From the experiment $\mathfrak{E}_{4,n,1}$ we construct an estimator
$\hat{\theta}_1$ for $\theta$. We define that
\[
\hat{\theta}_1    =   \sum_{k=1}^K m^{-1/2} \lambda_k^{-1} \int[\hat
{\Gamma}_1^{1/2} \varphi_k](t)\,dZ_1(t)   \varphi_k ,
\]
 where $K$ is an integer-valued smoothing parameter still to
be selected. Condition (\ref{eq:Cond_X}) guarantees that all $\lambda
_j$ are positive since, otherwise, $\lambda_j=0$ would yield that $E
|\langle X_1,\varphi_k \rangle|^2 = 0$ for all $k\geq j$, and hence
$\sum_{k\geq j} |\langle X_1,\varphi_k \rangle|^2 = 0$ a.s.\vspace*{-1pt} so that
$X_1$ would lie in the linear hull of $\varphi_1,\ldots ,\varphi_{j-1}$.
Thus the estimator $\hat{\theta}_1$ is well defined.

We introduce the data transformation $T_{4,5,n} \dvtx  \Omega_{4,n} \to
\Omega_{5,n}$ where
\begin{eqnarray*}
&&  T_{4,5,n}(\mathbf{x}_1,z_1,\mathbf{x}_2,z_2) \\
 &&  \qquad    =
\biggl(\mathbf{x}_1,z_1,\mathbf{x}_2,z_2 - (n-m)^{1/2} \int_0^\cdot[\hat{\Gamma
}_2^{1/2} \hat{\theta}_1](t)\,dt\\
&& \hspace*{62pt}
\qquad  \quad {} + (n-m)^{1/2} \int_0^\cdot[\Gamma
^{1/2} \hat{\theta}_1](t)\,dt \biggr) .
\end{eqnarray*}
 The transformation is fully accessible by the data drawn from
the experiment $\mathfrak{E}_{4,n,1}$ and the assumed knowledge of the
distribution of $\mathbf{X}$. We set $(\Omega_{5,n},\mathfrak{A}_{5,n}) =
(\Omega_{4,n},\mathfrak{A}_{4,n})$ where $\Omega_{4,n} = C_0^m([0,1])
\times C_0([0,1]) \times C_0^{n-m}([0,1]) \times C_0([0,1])$ and
$\mathfrak{A}_{4,n}$ is the corresponding Borel $\sigma$-algebra. The
data structure of $\mathfrak{E}_{4,n}$ is represented by $(\mathbf{X}_1,Z_1,\mathbf{X}_2,Z_2)$ when inserting the data set as an argument of
the mapping $T_{4,5,n}$. Note that $z_2=Z_2$ may be inserted in the
definition of the estimator $\hat{\theta}_1$. The integral occurring in
the definition of $\hat{\theta}_1$ is not defined for all continuous
functions $z_1$ but for almost all trajectories of $Z_1$. For the other
negligible trajectories the integral may conventionally be put equal to
zero to make the mapping $T_{4,5,n}$ well defined on the whole of its domain.

We define the experiment $\mathfrak{E}_{5,n}$ where ones observes the data
\[
(\mathbf{X}_1,Z_1,\mathbf{X}_2,Z_2')   =   T_{4,5,n}(\mathbf{X}_1,Z_1,\mathbf{X}_2,Z_2) ,
\]
 where the data $\mathbf{X}_1,Z_1,\mathbf{X}_2,Z_2$ are obtained
under experiment $\mathfrak{E}_{4,n}$. The experiment $\mathfrak
{E}_{5,n}$ is defined on the probability space $(\Omega_{5,n},\mathfrak
{A}_{5,n})$. Considering the definition of $T_{4,5,n}$, we realize that
the shift contained in the forth component is still available in the
experiment $\mathfrak{E}_{5,n}$ as the other components are kept.
Therefore, $T_{4,5,n}$ is an invertible transformation so that the
experiments $\mathfrak{E}_{4,n}$ and $\mathfrak{E}_{5,n}$ are equivalent.

In the experiment $\mathfrak{E}_{5,n}$, the component $Z_2'$ is still
an It\^{o} process conditionally on $\mathbf{X}_1,\mathbf{X}_2, Z_1$. Now we
introduce the experiment $\mathfrak{E}_{6,n}$ with $(\Omega
_{6,n},\mathfrak{A}_{6,n}) = (\Omega_{5,n},\mathfrak{A}_{5,n})$ where
one observes the data $\mathbf{X}_1, Z_1, \mathbf{X}_2$ and the It\^{o}
process $S_2(t)$, $t\in[0,1]$ with $S_2(0)=0$ and
\[
dS_2(t)   =   (n-m)^{1/2}  [\Gamma^{1/2} \theta ](t)\,dt   +
\sigma\, dW_2(t) .
\]
 In the notation of Section~\ref{S1.0}, we consider that
%
%e4.1 ###
\begin{eqnarray}\label{eq:4.1.t}
&& |E  R_{j,n,\theta}(\mathbf{X}_1,Z_1,\mathbf{X}_2,Z_2')   -   E
R_{j,n,\theta}(\mathbf{X}_1,Z_1,\mathbf{X}_2,S_2)  |\nonumber \\
&& \qquad    \leq  E\biggl | 1 - \exp\biggl ( - \sigma^{-1} \int\Delta_{5,6}(t)\,dW_2(t) - \frac1{2\sigma^2} \|\Delta_{5,6}\|_2^2 \biggr)
\biggr| \nonumber
\\[-8pt]
\\[-8pt]
 && \qquad
  \leq  2  E  \biggl\{1 - \exp \biggl(- \frac1{2\sigma^2} \|\Delta_{5,6}\|
_2^2 \biggr) \biggr\}^{1/2}\nonumber \\
  && \qquad    \leq  2  \biggl\{1 - \exp
 \biggl( - \frac1{2\sigma^2} E \|\Delta_{5,6}\|_2^2 \biggr) \biggr\}^{1/2} ,
\nonumber
\end{eqnarray}
 where
\begin{eqnarray*} \Delta_{5,6} &   = &  (n-m)^{1/2}  (\Gamma^{1/2}
\theta- \Gamma^{1/2} \hat{\theta}_1 + \hat{\Gamma}_2^{1/2} \hat{\theta
}_1 - \hat{\Gamma}_2^{1/2} \theta ) \\
 &   =&   (n-m)^{1/2} (\Gamma
^{1/2} - \hat{\Gamma}_2^{1/2}) (\theta- \hat{\theta}_1) .
\end{eqnarray*}
 Therein, we have used Girsanov's theorem, $\|R_{j,n,\theta}\|
_\infty\leq1$ for $R_{j,n,\theta} \in{\cal R}_{j,n,\theta}$ as
$j=5,6$, the Bretagnolle--Huber inequality and Jensen's inequality in
the last step.

Now we study the expectation occurring in (\ref{eq:4.1.t}) by
Parseval's identity with respect to the basis $\{\hat{\varphi}_{k,2}\}
_{k\geq1}$ and the orthogonal expansion of $\hat{\varphi}_{k,2}$ with
respect to $\{\varphi_j\}_{j\geq1}$ where $\{\hat{\varphi}_{k,2}\}
_{k\geq1}$ denotes the eigenfunctions of $\hat{\Gamma}_2$ and $\hat
{\lambda}_{k,2}$ the corresponding eigenvalues.
\begin{eqnarray*}
 E  \|\Delta_{5,6}\|_2^2   &=&   (n-m) \sum_{k=1}^\infty E  \Biggl|\sum
_{j=1}^\infty\langle\hat{\Gamma}_2^{1/2} \varphi_j - \Gamma^{1/2}
\varphi_j , \hat{\varphi}_{k,2}\rangle  \langle\varphi_j,\theta-\hat
{\theta}_1\rangle \Biggr|^2 \\
&     = &  (n-m) \sum_{k=1}^\infty E  \Biggl| \sum_{j=1}^\infty(\lambda
_j^{1/2} - \hat{\lambda}_{k,2}^{1/2}) \langle\varphi_j,\hat{\varphi
}_{k,2}\rangle  \langle\varphi_j , \theta- \hat{\theta}_1\rangle
\Biggr|^2 \\
&   \leq&  (n-m) \sum_{k=1}^\infty E \sum_{j=1}^\infty j^{-\gamma}
|\lambda_j - \hat{\lambda}_{k,2}| |\langle\varphi_j,\hat{\varphi
}_{k,2}\rangle|^2   \\
&&{}\times\sum_{j'=1}^\infty j'^\gamma E |\langle\varphi
_{j'} , \theta- \hat{\theta}_1\rangle|^2 ,
\end{eqnarray*}
 where we have used the Cauchy--Schwarz inequality for sums
and the elementary inequality $(\sqrt{x} - \sqrt{y})^2 \leq|x-y|$ for
all $x,y\geq0$. Therein $\gamma>0$ is still to be selected. Also, the
independence of $\hat{\Gamma}_2$ and $\hat{\theta}_1$ has been
utilized. Then, we apply the Cauchy--Schwarz inequality with respect to
the discrete random variable $V$ satisfying $P[V = |\hat{\lambda}_{k,2}
- \lambda_j|] = |\langle\varphi_j , \hat{\varphi}_{k,2}\rangle|^2$ for
all integers $k\geq1$ and some fixed integer $j$, conditionally on
$\mathbf{X}_2$. We conclude that
\begin{eqnarray*}
 E  \|\Delta_{5,6}\|_2^2   &\leq&  (n-m)  \Biggl \{\sum_{j'=1}^\infty
j'^\gamma E |\langle\varphi_{j'} , \theta- \hat{\theta}_1\rangle
|^2 \Biggr\} \\
 &&{}   \times E\sum_{j=1}^\infty j^{-\gamma}  \Biggl(
\sum_{k=1}^\infty|\lambda_j - \hat{\lambda}_{k,2}|^2  |\langle
\varphi_j,\hat{\varphi}_{k,2}\rangle |^2  \Biggr)^{1/2} \\
&   \leq&  (n-m)   \Biggl\{\sum_{j'=1}^\infty j'^\gamma E |\langle
\varphi_{j'} , \theta- \hat{\theta}_1\rangle|^2 \Biggr\}\\
&&{}\times \sum
_{j=1}^\infty j^{-\gamma}  \{E \|\hat{\Gamma}_2 \varphi_j - \Gamma
\varphi_j\|_2^2 \}^{1/2}  .
\end{eqnarray*}
We consider that
%
%e4.2 ###
\begin{eqnarray}\label{eq:4.0.m}
E\| \hat{\Gamma}_2 \varphi_j - \Gamma\varphi_j\|_2^2   &=&   E  \Biggl\|
\frac1{n-m} \sum_{k=1}^{n-m}  ( X_k \langle X_k,\varphi_j \rangle-
E X_k \langle X_k,\varphi_j \rangle ) \Biggr\|_2^2 \nonumber\\
 &
\leq&  (n-m)^{-1} E \|X_1\|_2^2 |\langle X_1,\varphi_j\rangle|^2 \nonumber\\
&   \leq&  (n-m)^{-1} c_j^2 \langle\Gamma\varphi_j , \varphi
_j\rangle  +   (n-m)^{-1} E \|X_1\|_2^4 1_{(c_j,\infty)}(\|X_1\|_2)
\\
&   \leq&  (n-m)^{-1} c_j^2 \lambda_j   +   (n-m)^{-1} \sum
_{k>c_j-1} (k+1)^4 P[\|X_1\|_2 \geq k] \nonumber\\
&   \leq & \mathrm{const.}\cdot(n-m)^{-1}  \bigl(c_j^2 \lambda_j   +
C_{X,0} \exp(-C_{X,1} c_j^{C_{X,2}} /2) \bigr)  \nonumber
\end{eqnarray}
 for $n$ sufficiently large where the sequence $(c_j)_j
\uparrow\infty$ remains to be determined. In order to obtain those
results, we impose the following:

\renewcommand{\thecond}{X}
\begin{cond}\label{condX} We assume that condition (\ref{eq:Cond_X})
holds true; $C_{\lambda,2} j^{-\alpha}   \geq  \lambda_j   \geq
C_{\lambda,1} j^{-\alpha}$ for all integer $j\geq1$ and some $\alpha
\geq2, C_{\lambda,2} > C_{\lambda,1} >0$; $EX_1 = 0$; $P[\|X_1\|_2
\geq x]   \leq  C_{X,0} \exp(-C_{X,1} x^{C_{X,2}})$ for all $x>0$
and some finite constants $C_{X,0},C_{X,1},C_{X,2}>0$.
\end{cond}

Condition~\ref{condX} imposes a polynomial lower bound on the sequence of the
eigenvalues of $\Gamma$. This assumption is very common in FLR (see,
e.g.,~\cite{CaiandHall2006}). When Condition~\ref{condX} is fixed the underlying
inverse problem can be viewed as a moderately ill-posed problem unlike
severely ill-posed problems where exponential decay of the eigenvalues
occurs. Condition~\ref{condX} also corresponds to the deconvolution setting with
ordinary smooth error densities in the related field of density
estimation based on contaminated data.

As an example for a stochastic process which satisfies Condition~\ref{condX}, we
mention the random variables
\[
X   =   \sum_{j=1}^\infty j^{-\alpha/2} G_j   \varphi_j ,
\]
 where the $\varphi_j$, integer $j$, form an arbitrary
orthonormal basis of $L_2([0,1])$; the $G_j$ are i.i.d. real-valued
centered random variables with a continuous distribution which is
concentrated on some compact interval, and $EG_1^2 = 1$. We stipulate
that $\alpha>2$. Easy calculations yield that the coefficients
$j^{-\alpha}$ and $\varphi_j$ are the eigenvalues and the eigenvectors
of the corresponding covariance operator $\Gamma$, respectively.
Stipulating that the sequence $\{\|\varphi_j\|_\infty\}_{j\geq1}$ is
bounded above (as satisfied, e.g., by the Fourier polynomials), we can
show that Condition~\ref{condX} is fulfilled. In particular, the random variable
$\langle g,X \rangle$ is continuously distributed for any $g \in
L_2([0,1])\setminus \{0\}$ since $\langle g,\varphi_j \rangle\neq0$ for at least one
integer $j$ so that the distribution of $\langle g,X \rangle$ is just
the convolution of an absolutely continuous distribution and some other
distribution; hence the distribution of $\langle g,\varphi_j \rangle$
has a Lebesgue density so that condition (\ref{eq:Cond_X}) can be
verified. All other assumptions contained in Condition~\ref{condX} can easily be
checked. Another even more important example for design distributions
are the Gaussian processes $X(t) = \int_0^t \sigma(s)\,dW(s)$, $t\in
[0,1]$, where $W$ denotes a standard Wiener process and $\sigma$ is a
sufficiently smooth function which is bounded from above and below by
positive constants. These processes satisfy Condition~\ref{condX} as well where
$\alpha=2$. The decay condition can be verified via the famous
reflection principle of Wiener processes.

Returning to the investigation of an upper bound on $E\|\Delta_{5,6}\|
_2^2$, the following inequality is evident:
%
%e4.3 ###
\begin{eqnarray} \label{eq:4.3.t}
E \|\Delta_{5,6}\|_2^2 &   \leq & (n-m)^{1/2}  \Biggl\{\sum
_{j'=1}^\infty j'^\gamma E |\langle\varphi_{j'} , \theta- \hat{\theta
}_1\rangle|^2 \Biggr\}\nonumber
\\[-8pt]
\\[-8pt]
 &&{}   \times \sum
_{j=1}^\infty j^{-\gamma}  \bigl(c_j^2 \lambda_j   +   C_{X,0} \exp
(-C_{X,1} c_j^{C_{X,2}}) \bigr)^{1/2} .
\nonumber
\end{eqnarray}
 We deduce that
%
%e4.4 ###
\begin{eqnarray} \label{eq:4.nneu}
 && E |\langle\varphi_{j'} , \theta- \hat{\theta}_1\rangle|^2\nonumber\\
 && \qquad    =
E \biggl|\langle\varphi_{j'} , \theta\rangle- 1_{\{j'\leq K\}} m^{-1/2}
\lambda_{j'}^{-1} \int[\hat{\Gamma}_1^{1/2} \varphi_{j'}](t)\,dZ_1(t) \biggr|^2 \nonumber\\
&& \qquad  = 1_{\{j'> K\}}  |\langle\varphi_{j'} , \theta\rangle |^2 \nonumber\\
&& \qquad  \quad   {}+ 1_{\{j'\leq K\}}\nonumber\\
&&\hphantom{ {}+} \qquad  \quad {}\times E  \biggl|\lambda_{j'}^{-1} \langle
\hat{\Gamma}_1 \theta,\varphi_{j'}\rangle- \langle\theta,\varphi_{j'}
\rangle+ \sigma m^{-1/2} \lambda_{j'}^{-1} \int[\hat{\Gamma}_1^{1/2}
\varphi_{j'}](t)\,dW_1(t) \biggr|^2 \nonumber\\
&& \qquad  = 1_{\{j'> K\}}  |\langle\varphi_{j'} , \theta\rangle |^2
\nonumber
\\[-8pt]
\\[-8pt]
& & \qquad  \quad  {}+ 1_{\{j'\leq K\}} \lambda_{j'}^{-2}  \{E  |
\langle\hat{\Gamma}_1 \theta,\varphi_{j'}\rangle- \langle\Gamma
\theta,\varphi_{j'} \rangle |^2 + \sigma^2 m^{-1} E \|\hat{\Gamma
}_1^{1/2} \varphi_{j'} \|_2^2 \}\nonumber \\
&& \qquad  \leq 1_{\{j'> K\}}  |\langle\varphi_{j'} , \theta\rangle
|^2\nonumber\\
&& \qquad  \quad {}+ 1_{\{j'\leq K\}} \lambda_{j'}^{-2}  \{m^{-1} \|\theta\|_2^2 E \|
X_1\|_2^2 |\langle X_1,\varphi_{j'}\rangle|^2 + \sigma^2 m^{-1} \lambda
_{j'} \} \nonumber\\
 && \qquad  \leq 1_{\{j'> K\}}  |\langle\varphi
_{j'} , \theta\rangle |^2 \nonumber\\
 && \qquad  \quad  {}+   1_{\{
j'\leq K\}}\nonumber\\
&&\hphantom{ {}+} \qquad  \quad {}\times \lambda_{j'}^{-2} \bigl \{m^{-1} \|\theta\|_2^2  \bigl(c_{j'}^2
\lambda_{j'}   +   C_{X,0} \exp(-C_{X,1} c_{j'}^{C_{X,2}}/2) \bigr) +
\sigma^2 m^{-1} \lambda_{j'} \bigr\} .
\nonumber
\end{eqnarray}

For further investigation of the asymptotic quality of the estimator
$\hat{\theta}_1$, some conditions on $P_X$ and $\Theta$ are required.
They are stated such that---combined with Condition~\ref{condX}---all previously
imposed assumptions concerning those characteristics are included.

\renewcommand{\thecond}{T}
\begin{cond}\label{condT} We assume that
\[
\sum_{k=1}^\infty(1+k^{2\beta})  |\langle\varphi_k , \theta\rangle
 |^2  \leq  C_\Theta
\]
 for all $\theta\in\Theta$ and some constants $\beta>(\alpha
+1)/2$ and $C_\Theta< \infty$, which are uniform with respect to
$\theta\in\Theta$.
\end{cond}

Condition~\ref{condT} says that the $\theta\in\Theta$ are uniformly well
approximable with respect to the orthonormal basis consisting of the
eigenfunctions of $\Gamma$. The parameter $\beta$ describes the degree
of this approximability. If the $\varphi_k$ were some Fourier
polynomials, then Condition~\ref{condT} could be interpreted as Sobolev
constraints on the set of the target functions.

We apply the parameter selection $K \asymp m^{1/(2\beta+\alpha+1)}$,
and we fix that $c_j = d_0 \log^{d_1}j$ with some constants $d_0, d_1$
sufficiently large and that $\gamma\in(0,\beta- \alpha/2 - 1/2)$.
Also, we choose\vadjust{\goodbreak} $m = \lfloor n/2 \rfloor$. Inserting that result into
(\ref{eq:4.3.t}), we deduce by Conditions~\ref{condX} and~\ref{condT} that
\[
\sup_{\theta\in\Theta} E \|\Delta_{5,6}\|_2^2   =   O \bigl(n^{(\alpha
+1-2\beta+2\gamma)/(4\beta+2\alpha+2)}\log^{d_3} n \bigr)   =   o(1)
\]
for some $d_3>0$, due to the inequality $\beta>(\alpha+1)/2$ and the suitable
selection of $\gamma$. Revisiting inequality (\ref{eq:4.1.t}), we have
finally proved by Section~\ref{S1.0}, paragraph (i) that the
experiments $\mathfrak{E}_{5,n}$ and $\mathfrak{E}_{6,n}$ are
asymptotically equivalent.

In the experiment $\mathfrak{E}_{6,n}$, the observation of $S_2$ allows
us to construct an estimator $\hat{\theta}_2$ for $\theta$ as well. It
is given by
\[
\hat{\theta}_2   =   \sum_{k=1}^K (n-m)^{-1/2} \lambda_k^{-1} \int
 [\Gamma^{1/2} \varphi_k ](t)\,dS_2(t) \varphi_k ,
\]
 where the parameter $K$ can be adopted from the estimator
$\hat{\theta}_1$. We specify the transformation $T_{6,7,n} \dvtx  \Omega
_{6,n} \to\Omega_{7,n}$ with
\begin{eqnarray*}
&& T_{6,7,n}(\mathbf{x}_1,z_1,\mathbf{x}_2,s_2) \\
 && \qquad    =
\biggl(\mathbf{x}_1,z_1 - m^{1/2} \int_0^\cdot[\hat{\Gamma}_1^{1/2} \hat{\theta
}_2](t)\,dt + m^{1/2} \int_0^\cdot[\Gamma^{1/2} \hat{\theta}_2](t)\,dt,\mathbf{x}_2,s_2 \biggr) .
\end{eqnarray*}
 Again the shift of the second component is accessible by the
other components which are maintained under the mapping so that
$T_{6,7,n}$ is invertible. Therefore, we define the experiment
$\mathfrak{E}_{7,n}$ by the observation of $T_{6,7,n}(\mathbf{X}_1,Z_1,\break\mathbf{X}_2,S_2)$
with $(\mathbf{X}_1,Z_1,\mathbf{X}_2,S_2)$ as under the experiment $\mathfrak{E}_{6,n}$. Hence, we put
$(\Omega_{7,n},\mathfrak{A}_{7,n}) = (\Omega_{6,n},\mathfrak{A}_{6,n})$
and obtain that $\mathfrak{E}_{7,n}$ is equivalent to $\mathfrak{E}_{6,n}$.

We define the experiment $\mathfrak{E}_{8,n}$ by the observation of
$(\mathbf{X}_1,S_1,\mathbf{X}_2,S_2)$ on the probability space $(\Omega
_{8,n},\mathfrak{A}_{8,n}) = (\Omega_{7,n},\mathfrak{A}_{7,n})$ where
$S_1$ denotes the It\^{o} process with $S_1(0)=0$ and
\[
dS_1(t)   =   m^{1/2}  [\Gamma^{1/2} \theta ](t)\,dt   +
\sigma\, dW_1(t) .
\]
 We can show that $\mathfrak{E}_{8,n}$ is asymptotically
equivalent to $\mathfrak{E}_{7,n}$ analogously to the proof of the
asymptotic equivalence of $\mathfrak{E}_{5,n}$ and $\mathfrak
{E}_{6,n}$. The only remarkable modification concerns the application
of the estimator $\hat{\theta}_2$ instead of $\hat{\theta}_1$. However,
even for that term we establish an upper bound at the same rate as for
estimator $\hat{\theta}_1$ since the asymptotic order of $m$ and $n-m$ coincide.

Taking a closer look at the data drawn from $\mathfrak{E}_{8,n}$, we
realize that the random variables $\mathbf{X}_1,S_1,\mathbf{X}_2,S_2$ are
independent. That occurs as we have replaced the empirical covariance
operators by the true deterministic one. Furthermore the data sets
$\mathbf{X}_1, \mathbf{X}_2$ do not carry any information on $\theta$ so that
$S_1,S_2$ represent a sufficient statistic for the whole empirical
information obtained under $\mathfrak{E}_{8,n}$. By Section~\ref{S1.0},
paragraph (ii), we conclude that $\mathfrak{E}_{8,n}$ is equivalent to
the experiment $\mathfrak{E}_{9,n}$ in which only the observations
$S_1,S_2$ are available. Thus we put $\Omega_{9,n} = C_0([0,1]) \times
C_0([0,1])$ and $\mathfrak{A}_{9,n}$ equal to the corresponding Borel
$\sigma$-algebra.\vadjust{\goodbreak}

We define the transformation $T_{9,10,n}\dvtx  \Omega_{9,n} \to\Omega
_{10,n}$ with $(\Omega_{10,n},\mathfrak{A}_{10,n}) = (\Omega
_{9,n},\mathfrak{A}_{9,n})$ by
\[
T_{9,10,n}(s_1,s_2)   =   \mathbf{A} (s_1,s_2)^T
\]
 with the matrix
\[
\mathbf{A}   =
\pmatrix{ m^{1/2}/n & (n-m)^{1/2}/n \vspace*{2pt}\cr
 m^{-1/2} & -(n-m)^{-1/2}
}
 .
\]
 We easily verify that $\mathbf{A}$ is invertible so that the
experiment $\mathfrak{E}_{10,n}$ which is defined by the observation of
$(T_1,T_2) = T_{9,10,n}(S_1,S_2)$ is equivalent to the experiment
$\mathfrak{E}_{9,n}$. We consider the characteristic function of the
$L_2([0,1])\times L_2([0,1])$-valued random variable $(T_1,T_2)$,
\begin{eqnarray*}
&& E \exp (i \langle t_1 , T_1 \rangle+ i \langle t_2 , T_2 \rangle
 ) \\
& & \qquad   =   E \exp (i \langle\mathbf{e}_1^T \mathbf{A}^T (t_1,t_2)^T , S_1
\rangle )  E\exp (i \langle\mathbf{e}_2^T \mathbf{A}^T (t_1,t_2)^T ,
S_2 \rangle ) \\
&& \qquad    =   \exp \biggl(i \biggl\langle t_1, \int_0^\cdot[\Gamma^{1/2} \theta](t)\,dt \biggr\rangle \biggr) \\
 && \quad \qquad {}\times\exp \biggl(-\frac1{2n} \sigma^2 \iint
t_1(x_1) \min\{x_1,x_2\} t_1(x_2)\,dx_1 \,dx_2 \biggr) \\
&& \quad \qquad {}\times\exp \biggl(-\frac1{2} \biggl[\frac1m + \frac1{n-m} \biggr] \sigma
^2 \iint t_2(x_1) \min\{x_1,x_2\} t_2(x_2)\,dx_1 \,dx_2 \biggr)
\end{eqnarray*}
 for any $t_1,t_2\in L_2([0,1])$ so that $T_1$ and $T_2$ are
two It\^{o} processes satisfying $T_1(0)=T_2(0)=0$ and
\begin{eqnarray*}
dT_1(t) &   = &  [\Gamma^{1/2} \theta](t)\,dt   +   n^{-1/2} \sigma\,
dW_3(t) , \\
dT_2(t) &   = &  \sigma \biggl(\frac1m + \frac1{n-m} \biggr)^{1/2}\,dW_4(t) ,
\end{eqnarray*}
 where $W_3$ and $W_4$ are two independent Wiener processes.
Thus $T_1$ and $T_2$ are independent, and $T_2$ is totally
uninformative with respect to the target function~$\theta$. Applying
Section~\ref{S1.0}, paragraph (ii) again, we have established
equivalence of $\mathfrak{E}_{10,n}$ and the experiment $\mathfrak
{E}_{11,n}$, which is equipped with $\Omega_{11,n} = C_0([0,1])$ and
the corresponding Borel $\sigma$-algebra $\mathfrak{A}_{11,n}$, and
characterized by the observation of the process $T_1$, which coincides
with the process $Y$ as defined in (\ref{eq:1.1.2}).

Summarizing we have shown asymptotic equivalence of the experiments
$\mathfrak{E}_{1,n}$ and $\mathfrak{E}_{11,n}$. That provides our final
main result, which will be given as a theorem below.

\begin{teo} \label{T:2}
Under the Conditions~\ref{condX} and~\ref{condT}, the FLR experiment $\mathfrak{E}_{1,n}$
with known design distribution and independent $N(0,\sigma
^2)$-distributed regression errors is asymptotically equivalent to the
white noise experiment $\mathfrak{E}_{11,n}$ where only the It\^{o}
process $Y$ as in (\ref{eq:1.1.2}) is observed.
\end{teo}

%s5 ###
\section{Sharp estimation for unknown $P_X$} \label{S:new}

We can combine our results with Theorem 1 in \cite
{CavalierandTsybakov2002}, which is due to~\cite{Pinsker1980}, in order
to derive a sharp minimax result with respect to the MISE for the FLR
problem under known covariance operator. It follows from there that
this sharp minimax risk corresponds to the sequence
\[
a_n   =   \sigma^2 n^{-1} \sum_{k=1}^\infty\lambda_k^{-1}  \bigl(1 -
\gamma(1+k^{2\beta})^{1/2} \bigr)_+ ,
\]
 where $\gamma$ is the unique solution of the equation
\[
\frac{\sigma^2}n \sum_{k=1}^\infty\lambda_k^{-1} (1+k^{2\beta})^{1/2}
\bigl (1 - \gamma(1+k^{2\beta})^{1/2} \bigr)_+   =   C_\Theta\gamma ,
\]
 under the conditions of Theorem~\ref{T:2}. More concretely,
there exists an estimator $\hat{\theta}$ of $\theta$ in the FLR model,
which satisfies
\[
\sup_{\theta\in\Theta} E \|\hat{\theta}-\theta\|_2^2   =   a_n
 \bigl(1+o(1) \bigr) .
\]
 Thus, any other estimator in the underlying model satisfies
the above equation when $=$ is replaced by $\geq$. We have established
sharp asymptotic constants.

Critically, we mention that the loss function $a_n^{-1} \|\hat{\theta}
- \theta\|_2^2$ is apparently not bounded. Still, asymptotic
equivalence yields coincidence of sharp minimaxity for the loss
function $\min\{D_n, a_n^{-1} \|\hat{\theta} - \theta\|_2^2\}$ for some
$(D_n)_n \to\infty$ sufficiently slowly. We can show that, in the
white noise inverse problem, the sharp constant result is extendable to
the truncated loss function. Using Theorem~\ref{T:2}, we have a sharp
lower bound for the FLR model even for the truncated loss function.

However, the design distribution $P_X$ is assumed to be known and
occurs in the definition of the minimax estimator. On the other hand,
the lower bound as derived from Theorem~\ref{T:2} in the previous
section provides a lower bound for the FLR model in the case of unknown
$P_X$ as well since nonknowledge of $P_X$ cannot improve this lower
bound. Thus if we succeed in showing that some estimator achieves these
asymptotic properties under the assumption of unknown $P_X$, then sharp
asymptotic minimaxity is extended to this more realistic condition.
Assuming that all conditions of Theorem~\ref{T:2} except the knowledge
of $P_X$ hold true, we propose the estimator
%
%e5.1 ###
\begin{equation} \label{eq:sh1}
\hat{\theta}   =   \sum_j w_j \frac1n \sum_{l=1}^n Y_l \langle
X_l,\hat{\varphi}_j\rangle\hat{\varphi}_j \hat{\lambda}_{j,\rho}^{-1}
\end{equation}
 for $\theta$ where $\hat{\lambda}_{j,\rho} = \max\{\hat
{\lambda}_j,n^{-\rho}\}$ for some $\rho\in(0,1/2)$. The weights $w_j$
remain to be specified. Using the techniques of the papers of \cite
{CaiandHall2006} and~\cite{HallandHorowitz2007}, the MISE of $\hat
{\theta}$ is equal to
%
%e5.2 ###
\begin{equation}\label{eq:sh2}
E \|\hat{\theta} - \theta\|_2^2    =
\sum_k E \biggl| w_k \frac{\hat{\lambda}_k}{\hat{\lambda}_{k,\rho}} - 1
\biggr|^2 |\langle\hat{\varphi}_k,\theta\rangle|^2   +   \frac{\sigma^2}n
\sum_k E w_k^2 \frac{\hat{\lambda}_k}{\hat{\lambda}_{k,\rho}^2} .
\end{equation}
 We stipulate that for $k> n^{\rho/\alpha} / \log n$ all
weights $w_k$ shall be put equal to zero. For all other $k$ we have
$\lambda_k \geq2 n^{-\rho}$ for $n$ sufficiently large so that
\begin{eqnarray*}
E \hat{\lambda}_k / \hat{\lambda}_{k,\rho}^2 &   \leq & \frac
1{(1-1/(2+\log k))\lambda_k}   +   n^{2\rho} \lambda_k   P [\hat
{\lambda}_k - \lambda_k < - \lambda_k /(2+ \log k )] \\[-2pt]
&   \leq&  \lambda_k^{-1} + \lambda_k^{-1}  \biggl(\frac1{\log k + 1} +
(2+\log k)^2 O(n^{2\rho-1}) \biggr) ,
\end{eqnarray*}
 where we have used Markov's inequality and that $E|\hat
{\lambda}_k - \lambda_k|^2$ is bounded from above by the expected
squared Hilbert--Schmidt norm of $\hat{\Gamma} - \Gamma$ and, hence, by
$O(n^{-1})$ (see, e.g.,~\cite{Bhatiaetal1983}). That requires the
following assumption:
%
%e5.3 ###
\begin{equation} \label{eq:as01}
\lambda_j - \lambda_{j+1}   \geq  \mathrm{const.}\cdot j^{-\alpha-1}
\end{equation}
(see also~\cite{HallandHorowitz2007}). We conclude that the
second term in equation (\ref{eq:sh2}) has the same asymptotic order as
\[
 \{1 + O(1/\log\log n) \}\cdot\frac{\sigma^2}n \sum_k w_k^2 \lambda
_k^{-1}  +   O (n^{-1} \log^{\alpha+1} n ) ,
\]
 under the above restriction with respect to the selection of
the weights. Focusing on the first term in (\ref{eq:sh2}), we deduce
by the Cauchy--Schwarz inequality that
\begin{eqnarray*}
\sum_k E \biggl| w_k \frac{\hat{\lambda}_k}{\hat{\lambda}_{k,\rho}} - 1
\biggr|^2 |\langle\hat{\varphi}_k,\theta\rangle|^2 &   \leq &  \biggl\{
\biggl(\sum_k E \biggl| w_k \frac{\hat{\lambda}_k}{\hat{\lambda}_{k,\rho}} -
1 \biggr|^2 |\langle\varphi_k,\theta\rangle|^2 \biggr)^{1/2} \\[-2pt]
&&\hphantom{ \biggl\{
\biggl(} {}   +   \mathrm{const.}\cdot \biggl(\sum_k E  |\langle\hat
{\varphi}_k - \varphi_k,\theta\rangle |^2 \biggr)^{1/2} \biggr\}^2 .
\end{eqnarray*}
 We consider that
\begin{eqnarray*}
  \sum_k E  |\langle\hat{\varphi}_k - \varphi_k,\theta\rangle
|^2   &=&   \sum_k E  \biggl|\sum_j \langle\hat{\varphi}_k - \varphi_k,
\varphi_j\rangle\langle\varphi_j,\theta\rangle \biggr|^2 \nonumber\\[-2pt]
&   \leq&  \mathrm{const.}\cdot C_\Theta\sum_{k,j} j^{-2\beta} E
|\langle\hat{\varphi}_k - \varphi_k,\varphi_j\rangle |^2 \nonumber\\[-2pt]
&   = &  \mathrm{const.}\cdot\sum_j j^{-2\beta} E |\langle\hat
{\varphi}_j - \varphi_j , \varphi_j\rangle |^2 \nonumber\\[-2pt]
&&{}  +   \mathrm
{const.}\cdot\sum_j j^{-2\beta} \sum_{k\neq j} E |\langle\hat
{\varphi}_k - \varphi_k , \varphi_j\rangle |^2 \\[-2pt]
&   \leq&  \mathrm{const.}\cdot\sum_j j^{-2\beta} E\|\hat{\varphi}_j -
\varphi_j \|_2^2 \nonumber\\[-2pt]
&&{}  +   \mathrm{const.}\cdot\sum_j j^{-2\beta} \sum
_{k\neq j} E |\langle\hat{\varphi}_k , \hat{\varphi}_j - \varphi
_j\rangle |^2 \nonumber\\[-2pt]
 \nonumber
&   \leq & \mathrm{const.}\cdot\sum_j j^{-2\beta} E\|\hat{\varphi}_j -
\varphi_j \|_2^2  \nonumber
\end{eqnarray*}
 by exploiting the orthonormality of the $\hat{\varphi}_j$ and
the $\varphi_j$ as well as Condition~\ref{condT} and Parseval's identity. Bhatia, Davis and McIntosh
~\cite{Bhatiaetal1983} provide that the squared
$L_2([0,1])$-distance between $\hat{\varphi}_j$ and $\varphi_j$ is
bounded from above by the squared Hilbert--Schmidt norm of $\hat{\Gamma
} - \Gamma$ multiplied by $8 j^{2\alpha+2}$ via condition~(\ref
{eq:as01}). Thus we have
%
%e5.4 ###
\begin{equation} \label{eq:sh3}
  \sum_k E  |\langle\hat{\varphi}_k - \varphi_k,\theta\rangle
|^2   =   O(n^{-1}) ,
\end{equation}
 where the constants contained in $O(\cdot)$ do not depend on
$\theta$ whenever
%
%e5.5 ###
\begin{equation} \label{eq:sh4}
\beta> \alpha+3/2 .
\end{equation}
 Returning to the consideration of the first term in (\ref
{eq:sh2}), we conclude that its asymptotic order reduces to that of
\[
\sum_k E \biggl|w_k \frac{\hat{\lambda}_k}{\hat{\lambda}_{k,\rho}} - 1
\biggr|^2 |\langle\varphi_k,\theta\rangle|^2 .
\]
 Then this term is bounded from above by
\begin{eqnarray*}
&&\sum_k  |w_k - 1 |^2 |\langle\varphi_k,\theta\rangle|^2   +
  \mathrm{const.}\cdot\sum_{k=1}^{\lfloor n^{\rho/\alpha}/\log n\rfloor
} |\langle\varphi_k,\theta\rangle|^2 \lambda_k^{-2}   E|\hat{\lambda
}_k - \lambda_k|^2 \\[-2pt]
 &&\hphantom{\sum_k  |w_k - 1 |^2 |\langle\varphi_k,\theta\rangle|^2   +
  \mathrm{const.}\cdot\sum_{k=1}^{\lfloor n^{\rho/\alpha}/\log n\rfloor
} }   {} +   O (n^{-2\beta\rho/\alpha
} (\log n)^{2\beta} ) \\[-2pt]
& & \qquad   \leq  O (n^{-2\beta\rho/\alpha} (\log n)^{2\beta} )  +
  \mathrm{const.}\cdot C_\Theta/n   +   \sum_k  |w_k - 1 |^2
|\langle\varphi_k,\theta\rangle|^2
\end{eqnarray*}
 by utilizing Condition~\ref{condT} and again the results of \cite
{Bhatiaetal1983}. The term $O (n^{-2\beta\rho/\alpha} )$ is
asymptotically negligible [i.e., bounded by $O(1/n)$] whenever $\rho>
\alpha/(2\alpha+3)$ as we have already imposed the condition (\ref
{eq:sh4}). It follows that the MISE of~(\ref{eq:sh1}) may be reduced to
its asymptotically efficient terms, that is,
%
%e5.6 ###
\begin{eqnarray} \label{eq:sh5}
E \|\hat{\theta} - \theta\|_2^2 &=& \{1 + o(1)\} \biggl(\sum_k |w_k-1|^2
|\langle\varphi_k,\theta\rangle|^2  + \frac{\sigma^2}n \sum_k w_k^2
\lambda_k^{-1} \biggr) \nonumber
\\[-9pt]
\\[-9pt]
&&{} + O (n^{-1} \log^{\alpha+1}n ) .
\nonumber
\end{eqnarray}
 The right-hand side of (\ref{eq:sh5}), however, corresponds
to the MISE of an oracle estimator which uses the\vadjust{\goodbreak} true versions of the
eigenvalues and eigenfunctions of $\Gamma$ instead of the empirical
ones. Also, it follows from~\cite{Pinsker1980} and \cite
{CavalierandTsybakov2002} that the estimator $\hat{\theta}$ as in (\ref
{eq:sh1}) attains the sharp asymptotic minimax risk when the weights
are chosen as
\[
w_k   =    (1 - \gamma\beta_k )_+ ,
\]
 when writing $\beta_k = (1+k^{2\beta})^{1/2}$ with an
appropriate deterministic parameter~$\gamma$. More precisely, we
consider $\gamma_n$ which we define by the unique zero of the function
$\Phi= \Phi_1 - \Phi_2$ with
\[
\Phi_1(x)    =   \sum_{k=1}^\infty\lambda_k^{-1} \beta_k  (1 - x
\beta_k )_+ , \qquad    \Phi_2(x)   =   C_\Theta x n / \sigma
^2
\]
 for $x\geq0$ where $\Phi_1$ and $\Phi_2$ are continuous
montonically decreasing and increasing, respectively. The selection
$\gamma= \gamma_n$ leads to asymptotic sharp optimality (see, e.g.,
\cite{CavalierandTsybakov2002}). Clearly, we have $\gamma_n \asymp
n^{-\beta/(2\beta+\alpha+1)}$. Otherwise, not even the convergence
rates are optimal as the required balance between the bias and the
variance term is violated. By condition (\ref{eq:sh4}) our additional
assumption saying that that $w_k\equiv0$ for $k> n^{\rho/\alpha} /
\log n$ is verified under this optimal selection of the weights when
stipulating that $\gamma> n^{-\beta/(2\beta+1)}$ as we have assumed
$\rho>\alpha/(2\alpha+3)$.

Still, the suggested selector is an oracle choice as it requires
knowledge of the true eigenvalues $\lambda_j$. That motivates us to
consider a data-driven selector $\hat{\gamma}$ of $\gamma$. First we
split the sample $(\mathbf{X},\mathbf{Y})$ into two independent data sets
$(\mathbf{X}_j,\mathbf{Y}_j)$, $j = 1,2$. The first data set $(\mathbf{X}_1,\mathbf{Y}_1)$
consists of $m$ pairs $(X_k,Y_k)$ where $m \asymp n(1 - 1/\log
n)$, and $(\mathbf{X}_2,\mathbf{Y}_2)$ contains all the other observations.
We employ $(\mathbf{X}_1,\mathbf{Y}_1)$ to estimate the function $\theta$
while the second data set (training data) is used to provide an
selector of $\gamma$. Concretely, we fix $\tilde{\gamma}$ as the unique
zero of $\hat{\Phi} = \hat{\Phi}_1 - \Phi_2$ where
%
%e5.7 ###
\begin{equation} \label{eq:weight} \hat{\Phi}_1(x)   =   \sum
_{k=1}^\infty (\hat{\lambda}'_{k,\rho} )^{-1} \beta_k  (1 - x
\beta_k )_+ .
\end{equation}
 Therein, $'$ indicates that the estimator is based on the
second data set. Then we define our selector of $\hat{\gamma}$ as $\mathrm
{med}\{n^{-\beta/(3\beta+1)},\tilde{\gamma},n^{-\beta/(2\beta+1)}\}$.
This truncation takes into account the a priori knowledge about the
true $\gamma_n$ so that $|\hat{\gamma}-\gamma_n|\leq|\tilde{\gamma
}-\gamma_n|$ almost surely for $n$ sufficiently large.

Thus determining $\hat{\gamma}$ does not require knowledge of $P_X$.
Now let us consider the MISE of the estimator $\hat{\theta}_{\hat{\gamma
}}$ where the index indicates the incorporated choice of the parameter
$\gamma$. By (\ref{eq:sh5}), we derive that
\begin{eqnarray*}
&&E \|\hat{\theta}_{\hat{\gamma}} - \theta\|_2^2\\
 && \qquad    =   o
\bigl(n^{-2\beta/(2\beta+\alpha+1)} \bigr)\\
&& \qquad  \quad {}  +   \{1+o(1)\} \\
 &&  \quad \qquad \hphantom{{}  +}{} \times
\Biggl(\sum_{k=1}^\infty|\langle\varphi_k,\theta\rangle|^2 E | (1 -
\hat{\gamma} \beta_k )_+ - 1 |^2 +   \frac{\sigma^2}m \sum
_{k=1}^\infty\lambda\lambda_k^{-1} E  (1 - \hat{\gamma} \beta_k
)_+^2 \Biggr) ,
\end{eqnarray*}
 where the terms contained in $o(1)$ do not depend on $\gamma
$. As the asymptotic order of $m$ and $n$ coincides, the estimator
based on $m$ data attains the same asymptotic rates and constants as
the estimator which uses even $n$ data, so our above calculations
remain valid. Therefore, the estimator $\hat{\theta}_{\hat{\gamma}}$
attains sharp minimax rates and constants whenever
%
%e5.8 ###
\begin{eqnarray}\label{eq:nnew1}
&&\sum_{k=1}^\infty E |(1-\hat{\gamma} \beta
_k)_+ - (1-\gamma_n \beta_k)_+ |^2 |\langle\varphi_k , \theta
\rangle|^2 \nonumber\\
 && \quad  {}  +   \frac1m \sum_{k=1}^\infty
\lambda_k^{-1} E |(1-\hat{\gamma} \beta_k)_+ - (1-\gamma_n \beta
_k)_+ |^2\\
  && \qquad    =   o \bigl(n^{-2\beta
/(2\beta+\alpha+1)} \bigr) ,\nonumber
\end{eqnarray}
 uniformly with respect to $\theta\in\Theta$. The first term
in (\ref{eq:nnew1}) is bounded from above by $\mathrm{const.}\cdot E|\hat
{\gamma}-\gamma_n|^2$. The second term has the upper bound
\[
O (c_n^{\alpha+2\beta+1} ) \cdot E|\hat{\gamma}-\gamma_n|^2   +
  O(1/n)\cdot\sum_{k>c_n n^{1/(2\beta+\alpha+1)}}^{\lceil\mathrm
{const.} \cdot n^{1/(2\beta+1)} \rceil} \lambda_k^{-1}   P[\hat{\gamma
}\leq\beta_k^{-1}]
\]
 for some sequence $(c_n)_n$ tending to infinity sufficiently
slowly. We deduce by Markov's inequality that (\ref{eq:nnew1}) is
satisfied if
\[
n^{\alpha/(2\beta+1)} \gamma_n^{-2\nu}\cdot E|\hat{\gamma}-\gamma
_n|^{2\nu}   +   E|\hat{\gamma}-\gamma_n|^2   =   o\bigl (n^{-2\beta
/(2\beta+\alpha+1)} \bigr)
\]
 for some fixed integer $\nu$. The assertion $|\tilde{\gamma
}-\gamma_n|>s_n$, for some positive-valued sequence $(s_n)_n \downarrow
0$ with $s_n / \gamma_n \to0$, implies that
\[
|\hat{\Phi}_1(\gamma_n+s_n) - \Phi_1(\gamma_n+s_n)| > C_\Theta s_n
n/\sigma^2
\]
or
\[
|\hat{\Phi}_1(\gamma_n-s_n) - \Phi_1(\gamma_n-s_n)| > C_\Theta s_n
n/\sigma^2 + |\Phi_1(\gamma_n) - \Phi_1(\gamma_n - s_n)| .
\]
 We have already imposed that $\rho> \alpha/(2\alpha+ 3)$ so
that $\lambda_k > n^{-\rho}$\vspace*{1pt} for all $k\leq(\gamma_n - s_n)^{-1/\beta
}$. That, however, yields
 $\|\hat{\Gamma}' - {\Gamma}\|_{\mathrm{HS}} \geq\mathrm{const.}\cdot
s_n  n^{-\rho} \gamma_n^{-1}$ where $\|\cdot\|_{\mathrm{HS}}$ denotes
the Hilbert--Schmidt norm of an operator.  Therein we have used the findings of \cite
{Bhatiaetal1983} again and the monotonicity of the functions $\hat{\Phi
}_1$,  $\Phi_1$, $\Phi_2$ as well as the definitions of
$\gamma_n$ and~$\tilde{\gamma}$. We deduce by Markov's inequality that
\begin{eqnarray*}
E|\hat{\gamma} - \gamma_n|^{2\nu} &   =&   s_n^{2\nu} + P [|\tilde
{\gamma}-\gamma_n|>s_n  ] \\
 &= &  s_n^{2\nu} + \mathrm{const.}\cdot n^{2\rho\mu} s_n^{-2\mu} \gamma
_n^{2\mu}   E \|\hat{\Gamma}' - {\Gamma} \|_{\mathrm{HS}}^{2\mu}
\end{eqnarray*}
 for any integer $\mu$. As all moments of $\|X_1\|_2$ are
finite by Condition~\ref{condT} we derive that
\[
E \|\hat{\Gamma}' - {\Gamma} \|_{\mathrm{HS}}^{2\mu}   =   O
\bigl((n-m)^{-\mu} \bigr) ,
\]
where we recall that $\hat{\Gamma}'$ is based on the training data set,
thus on $n-m \asymp n/\log n$ observations. As $\rho<1/2$ we conclude
by suitable choice of $(s_n)_n$ that
\[
E|\hat{\gamma} - \gamma_n|^{2\nu}   =   O ( [o_n \gamma_n
]^{2\nu} ) ,
\]
 where $(o_n)_n$ denotes some sequence tending to zero at an
algebraic rate. Choosing $\nu$ sufficiently large, we can finally
verify (\ref{eq:nnew1}) yielding the following proposition which
summarizes the investigation carried out in this section.
\begin{prop} \label{P:1}
We consider the FLR model in the setting of Theorem~\ref{T:2} except
the condition that $P_X$ is known. In addition, suppose that (\ref
{eq:as01}), (\ref{eq:sh4}) and $\rho\in (\alpha/(2\alpha+3),1/2
)$. Then, estimator (\ref{eq:sh1}) with the weight selector (\ref
{eq:weight}), which does not use $P_X$ in its definition, attains the
sharp minimax rate and constant with respect to the mean integrated
squared error; viewed uniformly over the function class $\Theta$ which
is defined via Condition~\ref{condT}.
\end{prop}

 Hence, under some additional conditions on the model, we have
established sharp minimaxity in the case where $P_X$ is unknown. Only
an arbitrary number between $\alpha/(2\alpha+3)$ and $1/2$ is supposed
to be known.

%s6 ###
\section{Discussion and conclusions} \label{S5}

We have proved equivalence of the FLR model and a white noise model
involving an empirical covariance operator in Theorem~\ref{T:1}. We
mention that $\sigma$ and $P_X$ can be treated as real nuisance
parameters in Section~\ref{S2}; more precisely, knowledge of those
quantities is not needed to apply the data transformations.

In contrast, for the asymptotic approximation in Section~\ref{S3},
$P_X$ must be known. Nevertheless, Section~\ref{S:new} shows that, with
respect to the MISE, the sharp asymptotic minimax risk can be taken
over to the case of unknown design distribution. Furthermore, under
specific parametric assumption on $P_X$, the condition of known $P_X$
can obviously be justified. Cai and Hall~\cite{CaiandHall2006}
explicitly mention Gaussian processes as examples for the random design
functions $X_j$. For instance, assuming that $X_j$ can be represented
as $X_j(t)  =   \int\xi(s)\,dW_j(s)$ with independent standard Wiener
processes $W_j$ as already suggested in the previous section, we
realize that the function $\xi$ is precisely reconstrucable based on
only one observation $X_1$. Then as $\xi$ is known the distribution
$P_X$ is known as well. Therefore, under this shape of $P_X$, the
assumption of known $P_X$ is not unrealistic at all. This phenomenon is
typical for the functional data approach and does not occur in
multivariate linear regression with finite-dimensional covariates. From
that point of view, the assumption of known design distribution causes\vadjust{\goodbreak}
less trouble in FLR compared to more standard regression problems.
Still this does not address the completely nonparametric case for $P_X$
and $\theta$.

As an interesting restriction, we have assumed that $\beta> (1+\alpha
)/2$ in Condition~\ref{condT}. Therefore, the quality of the approximation of the
target curve $\theta$ in the orthonormal basis consisting of the
eigenvalues of the covariance operator of the design variables must be
sufficiently high. If this basis consisted of Fourier polynomials then
that assumption could be interpreted as a smoothness condition on~$\theta$. That corresponds to the theorems in~\cite{Nussbaum1996} and
\cite{BrownandLow1996} where H\"{o}lder conditions are imposed, which
correspond to $\beta> 1/2$, in order to prove asymptotic equivalence
of the white noise model on one hand and density estimation and
nonparametric regression on the other hand. Otherwise, counterexamles
can be constructed (see~\cite{BrownandZhang1998}). To our best
knowledge our work represents the first proof of white noise
equivalence in a statistical inverse problem. It seems reasonable that
the essential condition is extended to $\beta> (1+\alpha)/2$ in this
setting as the selection $\alpha=0$ describes the setting of direct
estimation (noninverse problems). Still, the question of whether our
results are extendable to some $\beta\leq(1+\alpha)/2$ remains open.
In Section~\ref{S:new} we have studied the case of unknown $P_X$;
however, the regularity parameter $\beta$ is still assumed to be known.
Therefore, another interesting problem, which cannot be addressed
within the framework of this paper, is whether this sharp risk can be
achieved by an adaptive estimator, which does not use $\beta$ and
$C_\Theta$ in its construction. Approaches to adaptivity in FLR are
studied in~\cite{CaiandZhou2008}; however, that report seems to focus
on optimal rates rather than optimal constants.

Also, combining Theorem~\ref{T:2} and the results of Brown and Low \cite
{BrownandLow1996}, we conclude that, under reasonable conditions, the
FLR model is also equivalent to the standard nonparametric regression
problem, under which the data
\[
Y_j   =    [\Gamma^{1/2} \theta ](x_j)   +   \sigma
\varepsilon_j , \qquad  j=1,\ldots ,n ,
\]
 are observed where the $\varepsilon_j$ are i.i.d. and
$N(0,1)$-distributed, and the homogeneous fixed design setting $x_j =
j/n$, $j=1,\ldots ,n$, is applied.

\section*{Acknowledgments}
The author is grateful to Markus Reiss for a discussion on this paper
and to the reviewers for their inspiring comments.

%suskaldyti doi

% imsref loaded by smiklovaite, 2011-03-02 13:10:13
% imsref loaded by smiklovaite, 2011-03-02 13:14:33
%

\printaddresses

\end{document}